\documentclass[12pt]{amsart}
\usepackage{amsmath,amssymb}
\newtheorem{theorem}{Theorem}
\newtheorem{proposition}[theorem]{Proposition}
\title{Real and complex k-planes in convex hypersurfaces}

\author{Nikolai Nikolov}

\address{Institute of Mathematics and Informatics\\ Bulgarian Academy
of Sciences\\ Acad. G. Bonchev 8, 1113 Sofia,
Bulgaria}\email{nik@math.bas.bg}



\begin{document}

\begin{abstract} It is shown that that the rank of the second fundamental
form (resp.~the Levi form) of a $\mathcal C^2$-smooth convex
hypersurface $M$ in $\Bbb R^{n+1}$ (resp.~$\Bbb C^{n+1}$) does not
exceed an integer constant $k<n$ near a point $p\in M,$ then
through any point $q\in M$ near $p$ there exists a real
(resp.~complex) $(n-k)$-dimensional plane that locally lies on
$M.$
\end{abstract}

\maketitle

It is a classical result in the differential geometry that any
developable surface $M$ in $\Bbb R^3$ (i.e.~with zero Gaussian
curvature) is a part of a complete ruled surface (i.e.~through
every point of M there exists a straight line that lies on M).
Note that second fundamental form of such an $M$ has rank 0 or 1
at any point. A similar result holds in higher dimensions
(cf.~\cite[Lemma 2]{CL}):

(R) If the rank of the second fundamental form of a $\mathcal
C^2$-smooth hypersurface $M$ in $\Bbb R^{n+1}$ is a constant $k<n$
near a point $p\in M,$ then $M$ is locally generated by
$(n-k)$-dimensional planes. (In particular, if $k=0,$ then $M$ is
locally a hyperplane.)

This result has a complex version (see \cite[Theorem 6.1,
Corollary 5.2]{F}):

(C) If the rank of the Levi form of a $\mathcal C^2$-smooth real
hypersurface $M$ in $\Bbb C^{n+1}$ is a constant $k<n$ near a
point $p\in M,$ then $M$ is locally foliated by complex
$(n-k)$-dimensional manifolds. Moreover, if $k=0$ (i.e.~$M$ is
Levi-flat) and $M$ is real analytic, then $M$ is locally
biholomorphic to a complex hyperplane.

On the other hand, in both cases (real and complex), almost
nothing is known if the rank is not maximal and non-constant.

The aim of this note is consider the last case when the
hypersurface $M$ is convex, i.e. $M$ is a part of the boundary of
convex domain.

\begin{proposition}\label{main} The rank of the second fundamental
form (resp.~the Levi form) of a $\mathcal C^2$-smooth convex
hypersurface $M$ in $\Bbb R^{n+1}$ (resp.~$\Bbb C^{n+1}$) does not
exceed an integer constant $k<n$ near a point $p\in M$ if and only
if through any point $q\in M$ near $p$ there exists a real
(resp.~complex) $(n-k)$-dimensional plane that locally lies on
$M.$
\end{proposition}

\noindent\textsc{Remark.}~If $k=0$ in the complex case, then $M$
is locally linearly equivalent to the Cartesian product of $\Bbb
C^n$ and a planar domain (see \cite[Theorem 1]{BB}).

\begin{proof} If the respective real (complex) $(n-k)$-dimensional plane
exists for a point $q\in M,$ then the non-negativity of the second
fundamental form (the Levi form) at $q$ easily implies that the
rank of the form at $q$ does not exceed $k.$

For the converse, let first consider the complex case.

It is enough to show that through $p$ there exists a complex line
that locally lies on $M.$ Then, considering the intersection of
$M$ with the orthogonal complement of this line, we may proceed by
induction on $n$ to find $n-k$ orthogonal complex lines locally
lying on $M.$ The convexity of $M$ easily implies that the
$(n-k)$-dimensional planes, spanned by these lines, locally lies
on $M.$ Finally, note that the same holds for any point $q\in M$
near $p$ (since may replace $p$ by $q$).

Assume that there does not exist such a line. It is claimed in
\cite[p.~310]{S} and proved in \cite[Theorem 6]{NP} that $p$ is a
local holomorphic peak point for one of the sides, say $M^+$, of
$M$ near $p$ (the convex one). By \cite[Corollary 2]{B}, $p$ is a
limit of strictly pseudoconvex point of $M^+$ which is a
contradiction to the rank assumption.

The proof in the real case is similar. Recall that a point $q\in
M$ is called exposed if there exists a real hyperplane that
intersects $M$ in $p$ alone (i.e.~$p$ is a linear peak point). It
is enough to combine two facts:

$\bullet$ the set of exposed points is dense in the set of extreme
points (see \cite{St});

$\bullet$ the set of strictly convex points of $M^+$ (all the
eigenvalues of the second fundamental form are positive) is dense
in the set of exposed points.

The last fact can be shown following, for example, the proof of
\cite[Theorem]{B}.
\end{proof}

\noindent\textsc{Remark.} The 'only if' part ($\rightarrow$) of
Proposition \ref{main} remains true if we replace convexity by
real-analyticity. Indeed, if $c_M(q)$ denotes the rank of the Levi
form of $M$ at $q\in M$ and
$\displaystyle\widetilde{c}_M(p)=\limsup_{q\to p} c_M(q),$ then,
by (C), through any $q$ near $p$ with $c_M(q)=\widetilde{c}_M(p)$
there exists a complex line $(n-c_q)$-dimensional complex plane
that locally lies on $M.$ Then one may find a
$(n-c_q)$-dimensional complex plane with infinite order of contact
with $M$ at $p.$ Since $M$ is real-analytic, we conclude that this
plane lies on $M$ near $p.$ The real case follows analogously by
using (R) instead of (C).

\end{document}